\documentclass[12pt]{amsart}
\usepackage{amssymb,latexsym,amsmath}  

\newcommand{\K}{\mathcal{K}}
\renewcommand{\L}{\mathcal{L}}
\newcommand{\U}{\mathcal{U}}
\newcommand{\M}{{\mathcal{M}}}
\newcommand{\C}{\mathbb{C}}
\newcommand{\N}{\mathbb{N}}

\newcommand{\Td}{{\mathcal{T}_2}}

\newcommand{\ot}{\otimes}
\newcommand{\pf}{\noindent {\mbox{\textbf{Proof}.}} }
\newcommand{\ie}{\textit{i.e.\/}\ } 
\newcommand{\cst}{C$^*$} 

\newtheorem{thm}{Theorem}[section]

\newtheorem{lem}[thm]{Lemma}
\newtheorem{prop}[thm]{Proposition}

\theoremstyle{definition}
\newtheorem{defi}[thm]{Definition}

\newtheorem{rem}[thm]{Remark}

\newtheorem{ques}[thm]{Question}

\theoremstyle{remark}

\numberwithin{equation}{section}
\title{Deformations of infinite projections}
\author[Blanchard]{Etienne Blanchard}
\begin{document} 
\begin{abstract} 
Let $A=(A_x)$ be a (semi-)continuous field of \cst-algebras over a compact Hausdorff space $X$ and 
let $p=(p_x)$ be a projection in $A$ such that each $p_x\in A_x$ is properly infinite ($x\in X$). 
We prove that $p\oplus\ldots\oplus p$ ($l$ summands) is properly infinite in $M_l(A)$ for large enough $l\in\N$ if the $C(X)$-algebra $A$ is upper semi-continuous. 
But $p$ can be stably finite if $A$ is only lower semi-continuous.   
\end{abstract}

\maketitle
\section{Preliminaries}
A powerful tool in the classification of \cst-algebras is the study of their projections. 
\\ \indent 
Two projections $p, q$ in a \cst-algebra $A$ are said to be Murray-von Neumann \textit{equivalent} 
 (respectively (resp.) $p$ \textit{dominates} $q$) if there exists a partial isometry $v\in A$ with $v^*v=p$ and $vv^*=q$ 
 (resp. $v^*v\leq p$ and $vv^*=q$). 
 For short we write $p\sim q$ (resp. $q\preccurlyeq p$). 
The non-zero projection $p$ is said to be \textit{infinite} (resp. \textit{properly infinite}) if $p$ is equivalent to a proper subprojection $q< p$ (resp. $p$ is equivalent to two mutually orthogonal projections $p_1, p_2$ with $p_1+p_2\leq p$) and $p$ is \textit{finite} otherwise. 
\\ \indent 
J. Cuntz introduced the following generalization: A positive element $a$ in $A$ \textit{dominates} another positive element $b$ in $A$ (written $b\precsim a$) 
if and only if (iff) there exists a sequence $\{ d_n\}_n$ in $A$ such that $d_n^*ad_n^{}\to b$ (\cite{kiro}). 
Further $a\in A_+$ is called \textit{infinite} (resp. \textit{properly infinite}) iff 
there exists a non-zero positive element $b$ in $A$ such that $a\oplus b \precsim a\oplus 0$ in $M_2(A)$ (resp. $a\oplus a\precsim a\oplus 0$ in $M_2(A)$). 
And $a$ is said to be \textit{finite} if $a$ is not infinite. 
Kirchberg and R\o rdam proved that that these definitions coincide with the ones given in the previous paragraph in case $a$ is a projection (\cite[Lemma 3.1]{kiro}). 

Now a \cst-algebra $A$ is said to be infinite (resp. properly infinite) iff all strictly positive elements in A are infinite (resp. properly infinite). It is said to be finite (resp. stably finite) if all strictly positive elements in  $A$ are finite (resp. all strictly positive elements in $M_n(A)$ are finite for all positive integer~$n$).

\medskip 
In order to study deformations of such algebras, let us recall a few notions from the theory of $C(X)$-algebras. 
\\ \indent 
Let $X$ be a Hausdorff compact space and 
let $C(X)$ be the \cst-algebra of continuous functions on $X$ 
with values in the complex field $\C$. 
\begin{defi} 
A $C(X)$-algebra is a \cst-algebra $A$ endowed 
with a unital $*$--homo\-morphism from $C(X)$ 
to the centre of the multiplier \cst-algebra $\M(A)$ of~$A$. 
\end{defi}

For all $x\in X$, we denote by $C_x(X)$  the ideal 
of functions $f\in C(X)$ satisfying $f(x)=0$, 
by $A_x$ the quotient of $A$ 
by the \textit{closed} ideal $C_x(X) A$ and 
by $a_x$ the image of an element $a\in A$ in the \textit{fibre} $A_x$. 
Then the function 
\begin{equation} 
N(a):\,x\mapsto\| a_x\| =\inf\{\|\, [1-f+f(x)]a\|\,; f\in C(X)\}
\end{equation}
is upper semi-continuous by construction. 
The $C(X)$-algebra is said to be \textit{continuous} 
(or to be a \textit{continuous \cst-bundle over $X$}) if 
the function $x\mapsto\| a_x\|$ is actually continuous 
for all element $a$ in $A$.

\medskip 
\begin{defi} (\cite{bla2}) 
Given a continuous $C(X)$-algebra $B$, 
a \textit{$C(X)$-representation} of a $C(X)$-algebra $A$ on $B$ is a $C(X)$-linear map $\pi$ from $A$ 
to the multiplier \cst-algebra $\M(B)$ of $B$. 
Further $\pi$ is said to be a \textit{continuous field of faithful representations} if, 
for all $x\in X$, the induced representation $\pi_x$ of the fibre $A_x$ in $\M(B_x)$ is faithful. 
\end{defi} 
Note that the existence of such a continuous field of faithful 
representations $\pi$ implies that the $C(X)$-algebra $A$ is continuous 
since the function 
\begin{equation}\label{lsc}
x\mapsto \|\pi_x(a_x)\| =\sup\{\| (\pi (a)b)_x\|\,; 
b\in B \;\;\mbox{such that}\;\; \| b\|\leq 1\} 
\end{equation}
is lower semi-continuous for all $a\in A$. 
\\ \indent 
Conversely, any separable continuous $C(X)$-algebra $A$ admits 
a continuous field of faithful representations. 
More precisely, there always exists 
a unital positive $C(X)$-linear map $\varphi: A\to C(X)$ 
such that all the induced states $\varphi_x$ 
on the fibres $A_x$ are faithful (\cite{bla1}). 
By the Gel'fand-Naimark-Segal (GNS) construction this gives 
a continuous field of faithful representations of $A$ 
on the continuous \cst-bundle of compact operators $\K(E)$ 
on the Hilbert $C(X)$-module $E=L^2(A,\varphi)$.   

\bigskip 
A \textit{simple} \cst-algebra $A$ is purely infinite iff every non-zero hereditary \cst-subalgebra $B\subset A$ contains an infinite projection (\cite{Cu}). 
Possible generalisations to the non-simple case are the following: \\ 
-- A \cst-algebra $A$ is said to be \textit{purely infinite} (p.i.) 
iff $A$ has no non-zero character and 
for all $a, b\in A_+$, $\varepsilon>0$, with $b$ in the closed ideal of $A$ generated by $a$, 
there exists $d\in A$ with $\| b-d^* a d\|<\varepsilon$ (\cite{kiro}). 
\\
-- A \cst-algebra $A$ is said to be \textit{locally purely infinite} (l.p.i.) iff for all $b\in A$ and all ideal $J\triangleleft A$ with $b\not\in J$, there exists a stable \cst-subalgebra $D_J\subset \overline{b^* A b}$ such that $D_J\not\subset J$. 

Note that a \cst-algebra $A$ is p.i. 
iff for all $b\in A$, there exists a stable \cst-subalgebra $D\cong D\ot\K$ contained in the hereditary \cst-subalgebra $\overline{b^* A b}$ such that 
for all (closed two sided) ideal $J\triangleleft A$ with $b\not\in J$, then $D\not\subset J$ 
(\cite[prop. 5.4]{ro2}). 
Hence, every p.i. \cst-algebra is l.p.i. (\cite[prop. 4.11]{BKpi}). 
We shall study in this article a few problems linked to the converse implication. 
\\ \indent The author is grateful to E. Kirchberg and M. R\o rdam for helpful comments. 
He would also like to that the Humboldt University for invitations during which part of that work was written. 

\section{Continuous fields of properly infinite \cst-algebras} 
In this section, we study the stability properties of proper infiniteness under (upper semi-)continuous deformations. 
\\ \indent 
For all integer $n\geq 1$, $M_n(\C)$ is the \cst-algebra linearly generated by $n^2$ operators $\{ e_{i,j}\}$ satisfying the relations $e_{i,j} e_{k,l}=\delta_{j,k} e_{i,l}$ and 
\hbox{$(e_{i,j})^*=e_{j,i}$ ($1\leq i, j\leq n$)}. 
The Cuntz \cst-algebra $\mathcal{O}_n$ (resp. $\mathcal{T}_n$) is the unital \cst-algebra generated by $n$ isometries $s_1, \ldots, s_n$ satisfying the relation $s_1s_1^*+\ldots+s_ns_n^*=1$ \hbox{(resp. $s_1s_1^*+\ldots +s_ns_n^*\leq 1$).} \
\begin{defi} 
Given two \cst-algebra $A$ and $B$, a $*$-homomorphism $\pi: A\to B$ is said to be \textit{unit full} iff the closed two sided ideal generated by $\pi(A)$ in $B$ equals $B$. 
\end{defi} 

\begin{prop}\label{contpri} 
Let $X$ be a compact Hausdorff space and 
let $D$ be a unital separable $C(X)$-algebra the fibres of which are properly infinite.
Then $M_l(D)$ is properly infinite for some integer $l>0$. 
\end{prop} 

Let us first prove the following lemma which is essentially contained in \cite{BKpi}. 
\begin{lem}\label{lem2.2} 
a) Let $A,  B$ be unital \cst-algebras, $\pi: A\to B$ be a unital $*$-epimorphism, 
$\theta: \Td\to A$ and $\sigma:\Td\to B$ be unit full $*$-homomorphisms. 
Then there is a unit full $*$-homomorphism $\theta':\Td\to M_4(A)$ such that 
$(\imath\ot\pi)\theta'(r)=e_{1,1}\ot\sigma(r)$ for  all $r\in\Td$.

\noindent b) Suppose that the \cst-algebra $A$ is the pullback of the two unital 
\cst-algebras $A_1$ and $A_2$ along the $*$-epimorphisms $\pi_k: A_k\to B$ ($k=1, 2$). 
If $\theta_k:\Td\to A_k$ are unit full $*$-homomorphisms ($k=1, 2$), 
then there exists a unit full $*$-homomorphism $\widetilde{\theta}=(\widetilde{\theta}_1, \widetilde{\theta}_2): \Td\to M_4(A)\subset M_4(A_1)\oplus M_4(A_2)$ such that 
$\widetilde{\theta}_2(r)=e_{1,1}\ot\theta_2(r)$ for  all $r\in\Td$. 
\end{lem}
\pf a) Let $s_1, s_2$ be two isometries with orthogonal ranges generating the unital \cst-algebra $\Td$ and let $p\in\Td$ be the properly infinite projection $p=s_1 s_1^*+s_2 s_2^*$. 
Then the two full projections $\pi\theta(p)$ and $\sigma(p)$ are Murray-von Neumann equivalent in $B$ (\cite[lemma 4.15]{BKpi}). 
Thus, there exists a unitary $v\in M_2(B)$ with 
$$v^*\bigl(e_{1,1}\ot \pi\theta(p) \bigr)v=e_{1,1}\ot\sigma(p)\,.$$
Define the the unitary 
$u=1_{M_2(B)}-e_{1,1}\ot \sigma(p)\ + \sum_{k=1, 2} \bigl(e_{1,1}\ot\sigma(s_k)\bigr) v^*\bigl(e_{1,1}\ot\pi\theta(s_k^*) \bigr)v$ in $\U(M_2(B) )$. Then 
\begin{center} 
$e_{1,1}\ot\sigma(s_k)=u\, v^*\bigl(e_{1,1}\ot\pi\theta(s_k)\bigr)v$ \quad for $k=1, 2$. 
\end{center} 
Take unitary liftings $\widetilde{u}$ and $\widetilde{v}$ in $\mathcal{U}(M_4(A) )$ of the unitaries 
$u\oplus u^*$ and $v\oplus v^*$ which are in the connected component of the identity. 
The formulae $\theta'(s_k)=\widetilde{u}\,\widetilde{v}^*\bigl(e_{1,1}\ot\theta(s_k) \bigr)\widetilde{v}$ ($k=1, 2$) define a relevant $*$-homomorphism $\theta'$ from $\Td$ to $M_4(A)$. 

\noindent b)
The \cst-algebra $A$ is isomorphic to the \cst-subalgebra 
$\{ (a_1, a_2); a_j\in A_j \mathrm{\, and\,} \pi_1(a_1)=\pi_2(a_2)\}$ of $A_1\oplus A_2$. 
And $e_{1,1}\ot\pi_1\theta_1(p)=v^*(e_{1,1}\ot \pi_2\theta_2(p) )v$ for some $v\in\U(M_2(B) )$. 
Thus, by a), there exists an adequate $*$-morphism 
$\widetilde{\theta}=(\widetilde{\theta}_1, \widetilde{\theta}_2)$ from $\Td$ to $M_4(A)$. \qed 

\bigskip
\noindent\textbf{Proof of Proposition \ref{contpri}.} 
For all $x\in X$ there exist a open neighbourhood $U(x)$ of $x$ in $X$ with closure $\overline{U(x)}$ and a unital $*$-homomorphism $\Td\to D_{| \overline{U(x)}}:=D/C_0(X\setminus\overline{U(x)})D$ since $\Td$ is semiprojective (\cite[4.7]{blac}) and the fibre $A_x$ is properly infinite. 
Thus, there exist a finite covering $X=U_1\cup\ldots\cup U_n$ by open subsets $U_i$ and unital $*$-homo\-morphisms $\sigma_i: \Td\to D_{|\overline{U_i}}=:A_i$ ($1\leq i\leq n$). 
Now, step 2) of the above Lemma gives us a unit full $*$-homomorphism $\widetilde{\theta}: \Td \to M_l(D)$ for $l=4^{n-1}$, \ie such that 
the closed two sided generated by the projection  $q=\pi(1)$ equals $M_l(D)$. 

If we embed each $M_k(D)$ in $M_{k+1}(D)$ by $d\mapsto d\oplus 0$ ($k\in\N$), 
then $1_{M_l(D)}\preccurlyeq q\oplus\ldots\oplus q$ in $M_{n\,l}(D)$. 
But $q$ is properly infinite, i.e. $q\oplus q\preccurlyeq q$ (\cite{kiro}), and so 
$1_{M_l(D)}\oplus 1_{M_l(D)}\preccurlyeq q\oplus\ldots\oplus q\preccurlyeq q\leq 1_{M_l(D)}$.  \qed

\begin{rem} 
Uffe Haagerup indicated me another way to prove Proposition \ref{contpri}:  
If the unital \cst-algebra $D$ is stably finite \cst-algebra, 
then there exists a bounded non-zero lower semi-continuous quasi-trace on $D$ (\cite{bh}). 
Now, if $D$ is also a $C(X)$-algebra for some compact Hausdorff space $X$, 
this implies that there is a bounded non-zero lower semi-continuous quasitrace $D_x\to\C$ 
for at least some point $x\in X$ (\cite[Prop 3.7]{hrw}). 
But then, the fibre $D_x$ cannot be properly infinite. 
\end{rem}

\begin{ques} 
Does there exist a unital continuous $C(X)$-algebra $D$ the fibres of which are properly infinite  and which is finite? 
\end{ques}

\section{Lower semi-continuous fields of properly infinite \cst-algebras} 
Let us study whether the above results can be extended to lower semi-continuous (l.s.c.)  \cst-bundles $(A, \{\sigma_x\})$ over a compact Hausdorff space~$X$. 
Recall that any such \textit{separable} l.s.c. \cst-bundle admits a faithful $C(X)$-linear representation on a Hilbert $C(X)$-module $E$ such that, for all $x\in X$, the fibre $\sigma_x(A)$ is isomorphic to the induced image of $A$ in $\L(E_x)$ (\cite{lsc}). 
Thus, the problem boils down to the following: 
Given a separable Hilbert $C(X)$-module $E$ with infinite dimensional fibres $E_x$, 
the unit $p$ of the \cst-algebra $\L(E)$ of bounded adjointable $C(X)$-linear operators acting on $E$ has a properly infinite image in $\L(E_x)$ for all $x\in X$. 
But is the projection $p$ itself (properly) infinite in $\L(E)$? 
\\ \indent 
Dixmier and Douady have proved that this is always the case if the space $X$ has finite topological dimension (\cite{didou}). 
But it does not hold anymore in the infinite dimensional case: 
R\o rdam  has constructed an explicit example where $\L(E)$ is a finite \cst-algebra (\cite{ro}). 

\begin{ques} 
What happens if the compact Hausdorff space $X$ is contractible? 
\end{ques}

\noindent
\email{Etienne.Blanchard@math.jussieu.fr}\\
\address{IMJ,
Projet Alg\`ebres d'op\'erateurs, \indent
175, rue du Chevaleret, F--75013 Paris}

\end{document}